\documentclass[12pt]{article}

\usepackage[T2A]{fontenc}
\usepackage[utf8]{inputenc}

\usepackage{amsfonts}
\usepackage{amssymb}
\usepackage{amsmath}
\usepackage{amsthm}

\def \le {\leqslant}
\def \ge {\geqslant}

\theoremstyle{plain}

\begin{document}
 \begin{Large}
\centerline{\bf
Positive integers:
}
\centerline{\bf
counterexample to  W.M. Schmidt's conjecture
  }
 
 \vskip+1.5cm \centerline{ by {\bf Nikolay G. Moshchevitin}\footnote{ Research is supported by
the grant RFBR No. 09-01-00371-a}} \vskip+1.5cm
\end{Large}
\vskip+0.7cm
 
\begin{small} 
\centerline{\bf Abstract.}
We show that there exist real numbers $\alpha_1,\alpha_2$ linearly independent over $\mathbb{Z}$ together with 1
such that for every non-zero integer vector $(m_1,m_2)$  with $m_1\ge 0$ and $m_2\ge 0$
one has
$||m_1\alpha_1+m_2\alpha_2|| \ge 2^{-300} (\max(m_1, m_2))^{-\sigma}$ with
$\sigma   = 1.94696^+$.

\end{small}
\vskip+2.0cm

\section{Introduction}
Let $||\xi ||$ denotes the distance from real $\xi$ to the nearest integer.
Let
$\phi =\frac{1+\sqrt{5}}{2}
.$
In \cite{SCH}
W.M. Schmidt proved the following result.
\vskip+0.5cm
{\bf Theorem A.} (W.M. Schmidt)\,\,{\it
Let real numbers $\alpha_1,\alpha_2$
be linearly independent over $\mathbb{Z}$
together with 1. Then there exists a sequence of integer two-dimensional vectors
 $(x_1(i), x_2(i))$
 such that
}

1.\,\, $x_1(i), x_2(i) > 0$;

2.\,\, $||\alpha_1x_1(i)+\alpha_2x_2(i) ||\cdot (\max \{x_1(i),x_2(i)\})^\phi \to 0$ as $ i\to +\infty$.
\vskip+0.5cm

W.M. Schmidt
posed a conjecture that the exponent $\phi$ here may be  replaced by $2-\varepsilon$ with arbitrary positive $\varepsilon$
(see \cite{SCH1}).
In this paper we show this conjecture to be false.

Let $\sigma = 1.94696^+$ be the largest real root of the equation
\begin{equation}\label{sigma}
 x^4 - 2x^2-4x+1=0.
\end{equation}

{\bf Theorem 1.}\,\,
{\it There exist real numbers
 $\alpha_1,\alpha_2$
such that they are linearly independent over $\mathbb{Z}$
together with 1 and for every  integer vector
$(m_1,m_2)\in\mathbb{Z}^2$ with $m_1 , m_2 \ge 0$  and
$\max (m_1,m_2) \ge 2^{200}$
one has
$$
||m_1\alpha_1+m_2\alpha_2||\ge \frac{1}{2^{300}(\max (m_1,m_2))^\sigma}.
$$
}

We would like to formulate a related result from our paper \cite{moshe}.
For a real $\gamma \ge 2$ we define a function
$$
g(\gamma ) = \phi +\frac{2\phi - 2}{\phi^2\gamma -2}.
$$
One can see that $g(\gamma )$ is a strictly decresaing function and
$$
g(2) =2,\,\,\,\,\,\, \lim_{\gamma\to+\infty} g(\gamma ) = \phi.
$$
For positive $\Gamma$ define
$$
C(\Gamma ) = 2^{18} \Gamma^{\frac{\phi -\phi^2}{\phi^2\gamma-2}}.
$$

In \cite{moshe}   the following statement was proved. 
\vskip+0.5cm
{\bf Theorem B.} \,\,{\it
Suppose that real numbers $\alpha_1,\alpha_2$
satisfy the following Diophantine condition. 
For some $\Gamma\in (0,1)$ and $\gamma \ge 2$ the inequality
\begin{equation}
||\alpha_1m_1+\alpha_2m_2||\ge \frac{\Gamma}{(\max \{|m_1|,|m_2|\})^\gamma}
\label{bad}
\end{equation}
holds for all integer vectors $(m_1,m_2)\in \mathbb{Z}^2\setminus \{(0,0)\}$.
Then
 there exists an infinite sequence of integer two-dimensional vectors
 $(x_1(i), x_2(i))$
 such that}

1.\,\, $x_1(i), x_2(i) > 0$;

2.\,\, $||\alpha_1x_1(i)+\alpha_2x_2(i) ||\cdot (\max \{x_1(i),x_2(i)\})^{g(\gamma )} \le C(\Gamma )$ for all  $ i$.
\vskip+0.5cm
Of course   constants $2^{200},2^{300}$ and $2^{18}$ in Theorem 1 and in  the definition of $C(\Gamma )$ may be reduced.

\section{The construction}

We shall deal with the Euclidean norm for simplicity reason.
So we use $|\cdot |$ for the Euclidean norm of two- or three-dimensional vectors.
By ${\rm angle} ({\bf u}, {\bf v})$ we denote the angle between vectors ${\bf u}, {\bf v}$.

Define
\begin{equation}\label{tau}
 \tau = \frac{1+\sigma^2}{2\sigma} = 1.23029^+.
\end{equation}
Note that
\begin{equation}\label{tau1}
 \sigma\tau-1 >\tau.
\end{equation}
Put
\begin{equation}\label{omega}
 \omega = \tau + 1.
\end{equation}

{\bf Fundamental Lemma.}\,\,{\it
There exist
real numbers $\alpha_1, \alpha_2 \in \mathbb{R}$ 
linearly independent together with 1
over $\mathbb{Z}$ and such that there exisis a sequence
of integer vectors 
$$
{\bf m}_0 =(1,1,-1),\,\,\,\,\,
{\bf m}_\nu =(m_{0,\nu},m_{1,\nu}, m_{2,\nu}) \in \mathbb{Z}^3
, \nu =1,2,3,...
$$ satisfying the following conditions {\rm (i) -- (v)}.

{\rm (i)}\,\,
For any $\nu \ge 1$ the triple ${\bf m}_{\nu-1}, {\bf m}_\nu, {\bf m}_{\nu+1}$
 consists of linearly independent vectors, and each
two-dimensional sublattice
$${\cal L}_\nu = \langle 
{\bf m}_\nu, {\bf m}_{\nu+1}\rangle_\mathbb{Z}
$$
is complete,
that is
$${Z}^3 \cap
{
\rm span}\,
{\cal L}_\nu   = {\cal L}_\nu,\,\,\,\, \nu=0,1,2,3,... .                                     
$$

{\rm (ii)}\,\, 
Define
$$
\zeta_\nu  = m_{0,\nu} +m_{1,\nu}\alpha_1+m_{2,\nu}\alpha_2,\,\,\,\,\,
M_\nu = |\overline{\bf m}_\nu|.
$$
For every $\nu\ge 0$ one has
\begin{equation}\label{oo}
\frac{1}{2^5M_{\nu+1}^{\omega}}\le \zeta_\nu \le \frac{1}{M_{\nu+1}^{\omega}}.
\end{equation}

{\rm (iii)}\,\, 
$M_1\le 2^{100}$ and
for every $\nu\ge 1$ one has
\begin{equation}\label{q1}
 2^{10} M_\nu 
 \le M_{\nu+1}
\end{equation}
and 
\begin{equation}\label{q2}
H_\nu \le M_{\nu+1}
\le 2H_\nu,\,\,\,\ H_\nu= 
\frac{M_\nu^{\sigma\tau -1}}{2^{9}}.
\end{equation}

{\rm (iv)}\,\, For every $\nu\ge 0$ one has
$ m_{1,\nu}\cdot m_{2,\nu} <0$;
moreover for the vectors $${\bf e}_1 =\left(\begin{array}{c} 1 \cr 0\end{array}\right)
 ,\,\,\,\,\,       {\bf e}_1 =\left(\begin{array}{c} 0 \cr 1\end{array}\right)
 $$
and
$$
\overline{\bf m}_\nu =(m_{1,\nu},m_{2,\nu})\in \mathbb{Z}^2
$$
one has
\begin{equation}\label{eee}
{\rm angle} (\overline{\bf m}_\nu, \pm {\bf e}_j)  \ge \frac{1}{4},\,\,\,\,\,
j = 1,2.
\end{equation}

{\rm (v)}\,\, For every $\nu\ge 0$
for vectors
$$
\overline{\bf m}_\nu =(m_{1,\nu},m_{2,\nu}),\,\,\,\,
\overline{\bf m}_{\nu+1} =(m_{1,\nu+1},m_{2,\nu+1})
$$
 one has
$$
{\rm angle} (\overline{\bf m}_\nu, \pm \overline{\bf m}_{\nu+1})  \ge \frac{1}{4}.
$$

}

We give a sketched proof of Fundamental Lemma in
 Section 6.
It use standard argument related to
an inductive construction
of special singular (in the sense of A. Khintchine) vectors.
Inequality (\ref{tau1}) is of major importance.
Many different properties of singular vectors are discussed in our recent survey \cite{UMN}.

For every $\nu$ we define two-dimensional 
lattice $\Lambda_\nu = \langle 
\overline{\bf m}_\nu, \overline{\bf m}_{\nu+1}\rangle_\mathbb{Z}\subset \mathbb{Z}^2$.
Let $D_\nu$ be the fundamental volume of the lattice $\Lambda_\nu$.
Obviously
\begin{equation}\label{detup}
 D_\nu \le M_\nu M_{\nu+1}.
\end{equation}
From the condition {\rm (v)} one has
\begin{equation}\label{detlow}
 D_\nu \ge \frac{ M_\nu M_{\nu+1}}{2^5}.
\end{equation}

In the sequel we use the following notation.
  For an integer vector ${\bf m}=(m_0,m_1,m_2)\in\mathbb{Z}^3$ we
define
$$
\zeta = \zeta({\bf m})=m_0+m_1\alpha_1+m_2\alpha_2,
\,\,\,\,\,
\overline{\bf m} =
\overline {\bf m} ({\bf m}) = (m_1,m_2) \in \mathbb{Z}^2
$$
and
$$
M = M({\bf m}) =  |\overline{\bf m}|.
$$
 
In Sections 3,4,5 below we suppose that $\alpha_1, \alpha_2$
are the numbers from  Fundamental Lemma.

\section{Linearly independent vectors}

We prove a lemma concerning a lower bound for the value of $|\zeta ({\bf m})|$ 
in the case when the vector ${\bf m} \in \mathbb{Z}^3$ is linearly independent of 
vectors ${\bf m}_\nu, {\bf m}_{\nu+1}$.

Consider the segment
\begin{equation}\label{sege}
{\cal I}_\nu =
 \left[ (4M_\nu M_{\nu+1})^{1/\sigma}, M_{\nu+1}^\tau /8\right]
\end{equation}
(inequalities (\ref{tau1}) and (\ref{q1}) show that the left endpoitnt of the segment is
less than the right endpoint indeed).

{\bf Lemma 1.}
\,\,
{\it
Suppose that 
a vector ${\bf m} \in \mathbb{Z}^3$ is linearly independent of 
vectors ${\bf m}_\nu, {\bf m}_{\nu+1}$ and
\begin{equation}\label{i}
 M \in {\cal I}_\nu
.\end{equation}
Then
$$
|\zeta (m)| \ge M^{-\sigma}.
$$

}

Proof.

Consider the determinant
$$
\Delta 
=
\left|
\begin{array}{ccc}
m_{0} &m_{1} & m_{2}\cr
m_{0,\nu} &m_{1,\nu } & m_{2,\nu}\cr
m_{0,\nu+1} &m_{1,\nu +1} & m_{2,\nu+1}
\end{array}
\right|=
\left|
\begin{array}{ccc}
\zeta ({\bf m}) &m_{1} & m_{2}\cr
\zeta_{\nu} &m_{1,\nu } & m_{2,\nu}\cr
\zeta_{\nu+1} &m_{1,\nu +1} & m_{2,\nu+1}
\end{array}
\right|.
$$
We see from (\ref{oo}, \ref{omega}) that
$$
1\le |\Delta | \le 2|\zeta({\bf m})|M_\nu M_{\nu+1} + 4  MM_{\nu+1}^{-\tau}.
$$
From   the inequality $M \le M_{\nu+1}^\tau/8$ which follows from
(\ref{i}) we see that
$
 4 MM^{-\tau}_{\nu+1}\le 1/2.
$
That is why
$
|\zeta({\bf m})|M_\nu M_{\nu+1} \ge 1/4.
$ 
Now we take into account the lower bound for $M$ from (\ref{i}) and the lemma follows.$\Box$

\section{ Vectors dependent with $
{\bf m}_\nu, {\bf m}_{\nu+1}$}

Condition {\rm(i)} means  
 that   each integer vector ${\bf m} \in \mathbb{Z}^3$ which is linearly dependent together with
$
{\bf m}_\nu, {\bf m}_{\nu+1}$ can be written in a form
$$
{\bf m}=
\lambda {\bf m}_\nu +\mu {\bf m}_{\nu+1}$$
with integer $\lambda $ and $
\mu$.
So
if 
${\bf m} \in \mathbb{Z}^3$ is linearly dependent together with
$
{\bf m}_\nu, {\bf m}_{\nu+1}$ 
then for ``cutten'' vectors we have the equality
\begin{equation}\label{depp}
\overline{\bf m}=
\lambda \overline{\bf m}_\nu +\mu \overline{\bf m}_{\nu+1}
\end{equation}
with integer $\lambda $ and $
\mu$.

{\bf Lemma 2}.\,\,{\it
Suppose that the vector
${\bf  m} = (m_0,m_1,m_2)\in \mathbb{Z}^3$
satisfy the condition $m_1, m_2 \ge 0$
Suppose that 
vectors ${\bf m}, \bf{m}_\nu, {\bf m}_{\nu+1}$ are linearly dependent for some $\nu$.
Then
$$
|\zeta ({\bf m})| \ge 2^{-300}M^{-\sigma}.
$$}

Proof.

We
can
split
two-dimensional lattice $\Lambda_\nu$ 
into a countable union of one-dimensional lattices 
$\Lambda_{\nu,\mu}$ in the following way:
$$
\Lambda_\nu 
=\bigsqcup_{\mu\in \mathbb{Z}} \Lambda_{\nu,\mu},\,\,\,
\Lambda_{\nu, \mu} =\{ {\bf z} =(z_1,z_2)\in \Lambda_\nu:\,\,{\bf z} =\lambda \overline{\bf m}_\nu+\mu
\overline{\bf m}_{\nu+1},\,\, \lambda\in \mathbb{Z}\}.
$$
By the condition {\rm (iv)} there
is no non-zero points  $ (z_1,z_2)\in \Lambda_{\nu,0}$ satisfying $z_1\cdot z_2 \ge 0$.

Suppose that $\mu\neq 0$.
As the fundamental volume of $
\Lambda_\nu$ is equal to $D_\nu$ we see that the Euclidean distance between any two neighbouring
lines ${\rm aff}\, \Lambda_{\nu,\mu}$ and
${\rm aff}\, \Lambda_{\nu,\mu+1}$
is equal to $ D_\nu /\sqrt{m_{1,\nu}^2 +m_{2,\nu}^2}$.
That is why the conditions
\begin{equation}\label{o}
(m_1,m_2) \in  \Lambda_{\nu,\mu},\,\,\, m_1, m_2 \ge 0
\end{equation}
imply
\begin{equation}\label{o1}
\max (m_1,m_2) \ge \frac{|\mu|D_\nu}{2M_\nu} \ge \frac{|\mu|M_{\nu+1}}{2^8}
\end{equation}
(in the last inequality we use (\ref{detlow})).
 
From the other hand conditions (\ref{o}) together with (\ref{eee}) from {\rm (iv) }
lead to the inequality 
$$|\lambda |\ge \frac{|\mu|M_{\nu+1}}{2^8M_\mu}$$
 for the coefficient $\lambda$
from (\ref{depp}).
So (we apply lower bound from {\rm (ii)} for $\zeta_\nu$ and  upper bound from {\rm (ii)} for 
$\zeta_{\nu+1}$) we see that
$$
|\zeta ({\bf m})| =|\lambda\zeta_\nu+\mu\zeta_{\nu+1}|\ge
|\mu| \left(\frac{ M_{\nu+1}}{2^8M_\nu}\zeta_\nu -\zeta_{\nu+1}\right)
\ge |\mu|\left(\frac{1}{2^{13}M_\nu M_{\nu+1}^\tau }-\frac{1}{M_{\nu+2}^{\tau+1}}\right)
.
$$
 We apply lower bound (\ref{q1}) from {\rm (iii)} to get
 $$|\zeta({\bf m})| \ge \frac{|\mu|}{2^{14} M_\nu M_{\nu+1}^\tau}.
$$
Now lower bound (\ref{q2}) from {\rm (iii)} gives
\begin{equation}\label{o3}
 |\zeta({\bf m})| \ge2^{-14-\frac{100}{\sigma\tau -1}} |\mu| M_{\nu+1}^{\tau+\frac{1}{\sigma\tau -1}}
=2^{-14-\frac{100}{\sigma\tau -1}} |\mu |M_{\nu+1}^{-\sigma}>2^{-200}|\mu |M_{\nu+1}^{-\sigma}
\end{equation}
(here we use the definition of $\sigma$ as a root of (\ref{sigma}) and  (\ref{tau}) to see that $\tau + \frac{1}{\sigma\tau -1} =\sigma$).
Now we  combine (\ref{o1},\ref{o3}) and (\ref{q2}) from  the condition {\rm (iii})
to get
$$
|\zeta({\bf m})| \ge \frac{|\mu|^{1+\sigma}}{2^{300}M^\sigma} >
\frac{1}{2^{300}M^\sigma}.
$$
Lemma 2 is proved.$\Box$

\section{Proof of Theorem 1}

We take $\alpha_1, \alpha_2$ from Fundamental Lemma. Consider an integer vector ${\bf m} = (m_0,m_1,m_2)$ with $m_1,m_2\ge 0$.
We may suppose that $|\zeta({\bf m})| =||m_1\alpha_1+m_2\alpha_2||.$
If for some $\nu$ vectors 
\begin{equation}\label{bf}{\bf m}, {\bf m}_\nu,{\bf m}_{\nu+1}
\end{equation}
 are linearly dependent then
application of Lemma 2 proves Theorem 1.
So we may suppose that all
triples (\ref{bf}) consist of linearly independent vectors for every $\nu \ge 0$.
Now to prove Theorem 1 we may use Lemma 1.
It is enough to show that
$$
\bigcup_{\nu \ge 0} {\cal I}_\nu \supset [2^{200}, +\infty)
$$
(segments ${\cal I}_\nu$ are defined in (\ref{sege})).
But this follows from the condition $ M_1 \le 2^{100}$ and the inequality
$$
  (4M_\nu M_{\nu+1})^{1/\sigma}\le M_{\nu}^\tau /8
.$$
The last inequality is a corollary of the right inequality from (\ref{q2}).$\Box$

\section{ Fundamental Lemma: sketch of a proof }

Let ${\bf m}\in \mathbb{Z}^3$ be an integer vector.
The formulation of Fundamental Lemma deals 
with the values of $M= M({\bf m}) = |\overline{\bf m}|$.
To describe the ideas of the proof it is much more convenient to consider
the Euclidean norm
${\rm M} = |{\bf m}|$ of the vector ${\bf m}$ itself than the Euclidean norm 
$M = |\overline{\bf m}|$  of the ``cutten''
vector $\overline{\bf m}\in \mathbb{Z}^2$.
Of cource values of $M$ and ${\rm M}$ are of the same order for all integer vectors ${\bf m}$
under consideration.
We 
  may assume that
$M\le {\rm M}\le 2M$.

Let
$$
\hbox{\got S} =\{{\bf x}=(x_0,x_1,x_2)\in \mathbb{R}^3:\,\, |{\bf x}|=1\}
$$
be the unit sphere.
We construct a sequence of nested closed sets  ${\cal B}_\nu \subset \hbox{\got S} $ by induction.
Their unique common point $
{\bf x}^*=(x_0^*,x_1^*,x_2^*) \in \bigcap_\nu {\cal B}_\nu$ will define 
real  numbers $\alpha_1=x_1^*/x_0^*, \alpha_2 = x_2^*/x_0^*$ which
  satisfy the conclusion of Fundamental Lemma.

The base of inductive process is trivial.

To proceed the inductive step we suppose that the following objects   are alredy constructed:
\vskip+0.4cm
1) primitive integer vectors ${\bf m}_j= (m_{0,j}, m_{1,j}, m_{2,j}) ,\,\,\, 0\le j \le \nu$
with ${\rm M}_j = |{\bf m}_j|$; we suppose
that these vectors  satisfy conditions {\rm (iv), (v)};

2)   vectors  ${\bf \xi}_j =(\xi_{0,j},\xi_{1,j},\xi_{2,j})\in \hbox{\got S}$ such that
${\bf \xi}_j \perp {\bf m}_j,\,\,\, {\bf \xi}_j \perp {\bf m}_{j+1},\,\,\,0\le j \le \nu-1$
,
and one-dimensional linear subspases $\Xi_j = {\rm span }\, \xi_j,\,\,\, 0\le j \le \nu -1$
;

3) 
two-dimensional linear subspaces
$$
\ell_j^0 = \left\{ (x_0,x_1,x_2) \in \mathbb{R}^3 :
\, m_{0,j}x_0+m_{1,j}x_1+ m_{2,j} x_2 =0\right\},\,\,\,\
0\le j \le \nu-1,
$$
and
two-dimensional  affine subspaces 
$$
\ell_j^1 = \left\{ (x_0, x_1,x_2) \in \mathbb{R}^3 :
\, m_{0,j}x_0+m_{1,j}x_1+ m_{2,j} x_2 =\frac{1}{2{\rm M}_{j+1}^\omega}\right\},\,\,\,\
0\le j \le \nu-1;
$$


4) 
cylinders
$$
{\cal C}_j =\left\{
{\bf x} \in \mathbb{R}^3:\,\,\, {\rm dist}\, ({\bf x}, \Xi_j)\le  \frac{1}{{\rm M}_{j+1}^{\omega}{\rm M}_j}
\right\},\,\,\, 0\le j \le \nu-1                                                                                          
$$
(here ${\rm dist}\, (\cdot,\cdot )$ denotes the Euclidean distance between sets) and closed
sets
$$
{\cal G}_j =  
\left\{ {\bf x} \subset \hbox{\got S} \cap {\cal C}_j:\,\,
m_{0,j}x_0+m_{1,j}x_1+ m_{2,j} x_2 \ge \frac{1}{2{\rm M}_{j+1}^\omega }
 \right\}
\subset \hbox{\got S},\,\,\, 0\le j \le \nu-1
$$
(so a part of the boundary of ${\cal G}_j $ belongs to $\ell_j^1$);

5) two-dimensional 
complete sublattices ${\cal L}_j = \langle {\bf m}_j, {\bf m}_{j+1}\rangle_\mathbb{Z}, \,\,\, j =0,..., \nu-1$
with fundamental volumes $d_j$ satisfying inequalities
\begin{equation}\label{volume}
\frac{{\rm M}_{j}{\rm M}_{j+1}}{2^{5}}\le d_j \le {\rm M}_{j}{\rm M}_{j+1}
\end{equation}
(here the right inequality  is trivial, the left one means that the angle between
vectors ${\bf m}_{j-1}, {\bf m}_{j}$ is bounded from below);

6) we suppose that the vector ${\bf m}_\nu $ is defined,
so we can consider linear subspace
$$
\ell_\nu^0 = \left\{ (x_0,x_1,x_2) \in \mathbb{R}^3 :
\, m_{0,\nu}x_0+m_{1,\nu}x_1+ m_{2,\nu} x_2 =0\right\};
$$
we suppose that linear subspases $\ell_j^0$ for every $j$ from the range
$1\le j \le \nu$
satisfy the condition
$$
\ell_j^0 \cap {\cal G}_{j-1}
\neq \varnothing
,\,\,\, 1\le j\le \nu;
$$
moreover we suppose that for any $j$ from the range $ 1\le j \le \nu$
there is a point 
${\bf \eta}_j =(\eta_{0,j},\eta_{1,j}.\eta_{2,j})\in \ell_j^0\cap {\cal G}_{j -1}$
such that  the set
$$
{\cal B}_j =\left\{ {\bf x} \in \hbox{\got S}:\,\,\, 
|{\bf x}- {\bf  \eta} \,  _j|
\le 
\frac{1}{2^6{\rm M}_j^\omega {\rm M}_{j-1}}
\right\}
$$
satisfy the condition
\begin{equation}\label{nested1}
{\cal B}_j \subset {\cal G}_{j-1} \subset {\cal B}_{j-1}.
\end{equation}
 Here we should note that
 ${\cal B}_0\supset {\cal B}_1 \supset \cdots \supset {\cal B}_{\nu-1} \supset {\cal B}_\nu$.

\vskip+0.4cm
We suppose that vectors ${\bf m}_j ,\,\, 1\le j \le \nu $
and every  couple $\alpha_1, \alpha_2$ of the form $ \alpha_1 = x_1/x_0, \alpha_2=x_2/x_0,\,\,
{\bf x} = (x_0, x_1,x_2)\in B_{\nu-1}$
satisfy all the conditions {\rm(i) -- (v)} of Fundamental Lemma
which are defined
up to the $(\nu -1)$-th step.

Our task is to define an integer vector
${\bf m}_{\nu+1}$
and all related objects of the $\nu$-th step.

Consider ${\bf n} =(n_0,n_1,n_2) \in \mathbb{Z}^3$  such that the triple
${\bf n}, {\bf m}_{\nu-1}, {\bf m}_\nu$ form a basis of $\mathbb{Z}^3$.
Such vector does exist as the lattice ${\cal L}_{\nu-1}$ is complete.
We may suppose that
 \begin{equation}\label{enny}
\max (|n_1|,|n_2| )\le  {\rm M}_\nu.
\end{equation}
We  consider two-dimensional lattices
$$
{\cal L}_{\nu-1,\mu}
=\{
{\bf z} = \lambda_1{\bf m}_{\nu-1} +\lambda_1 {\bf m}_\nu + \mu {\bf n},
\,\,\,\, \lambda_1,\lambda_2 \in \mathbb{Z}\}.
$$
Note that
$$\mathbb{Z}^3 = \bigsqcup_{\mu\in \mathbb{Z}} {\cal L}_{\nu-1,\mu}.
$$
In fact $
{\bf n } \in {\cal L}_{\nu-1,1}
$.
The Euclidean distance between the neighbouring affine subspaces
${\rm aff}\, {\cal L}_{\nu-1,\mu}$ and 
${\rm aff}\, {\cal L}_{\nu-1,\mu+1}$
is equal to $d_{\nu-1}^{-1}$.
Put
$$
\mu_* = d_{\nu-1} H_\nu {\rm M}_{\nu}^{-\omega}{\rm M}_{\nu-1}^{-1}
$$
In fact $\mu_0$ is of  the size
$$
  \mu_* \asymp M_\nu^t, \,\,\,\, t = \sigma \tau - \omega>0 
$$
(here the last inequality follows from (\ref{tau1})).

Now  here we define two-dimensional linear subspace
$
\ell_\nu^*
\subset \mathbb{R}^3$
and a point ${\bf w}_\nu \in {\rm aff}\, {\cal L}_{\nu-1,\mu_*}$ by the following way.
Consider the unique one-dimensional affine subspace $\pi \subset \mathbb{R}^3$
such that 

1) $\eta_\nu \in \pi$,

2) $\pi$ is parellel to $\ell_{\nu -1}^0$,

3) the intersection $\pi \cap \hbox{\got S}$ consists of just one point 
$\eta_\nu$.

We define $\ell_\nu^* $ as follows:
$
\ell_\nu^* = {\rm span}\, \pi.
$
Let ${\bf w}_\nu \in {\rm aff} \, {\cal L}_{\nu-1,\mu_*} $ be the unique point such that
$ {\bf w}_\nu \perp \ell_\nu^*$.
Now we define the disk
$$
{\cal D}_\nu =\left\{ {\bf w}\in {\rm aff}\,{\cal L}_{\nu-1,\mu_*}:\,\,\,
|{\bf w}- {\bf w}_\nu|\le \frac{H_\nu}{2^{10}}\right\}.
$$
Easy calculation shows that
$$
{\bf w} \in {\cal D}_\nu \,\,\Longrightarrow\,\,
H_\nu \le |\overline{\bf w}|\le 2 H_\nu.
$$
For $ {\bf w}=(w_0,w_1,w_2)$ 
we consider two-dimensional linear subspase
$$
\ell[{\bf w}] = \{
{\bf x} = (x_0,x_1,x_2)\in \mathbb{R}^3:\,\,\,
w_0x_0+w_1x_1+w_2x_2 = 0\}.
$$
Consider a smaller ball ${\cal B}_\nu' \subset {\cal  B}_\nu$ with the same center ${\bf \eta}_\nu$
 and radius $ \frac{1}{2^7{\rm M}_\nu^{-\omega} {\rm M}_{\nu-1}}$. 
Easy calculation shows that 
 $$
 {\bf w}\in {\cal D}_{\nu}
\,\,\Longrightarrow \,\, \ell [{\bf w}]\cap \ell_\nu^0\cap {\cal B}_\nu' \neq \varnothing.
$$

Now if we   take an integer vector ${\bf m}_{\nu+1} \in {\cal D}_{\nu}$
 the conditions {\rm (ii), (iii)}  are satisfied for $\nu$-th step.

Vector ${\bf \xi}_{\nu}$, subspaces $\ell_\nu^1, \ell_{\nu+1}^0$ and the set  ${\cal G}_{\nu}$
 are defined automatically.
We can easily take ${\bf \eta}_{\nu+1}$ with all necessary properties,
in particular  we condtruct ${\cal B}_{\nu+1}$
(the second embedding  in (\ref{nested1}) with $j =\nu+1$ follows from the largeness of the value of $M_{\nu+1}$,
the first one can be ensured as the angle between vectors ${\bf m}_{\nu+1}. {\bf m}_{\nu}$
is almost the same as the angle between vectors    ${\bf m}_{\nu}. {\bf m}_{\nu-1}$).

Now we must explain how to ensure the condition
{\rm (i)}. and {\rm (iv), (v)}.

To get {\rm (i)} we should note that a vector ${\bf n} = (n_0,n_1,n_2)\in {\cal L}_{\nu-1,1}$ which
completes the pair ${\bf m}_{\nu-1}, {\bf m}_\nu$ to a basis of $\mathbb{Z}^3$
may be found in any box of the form
$$
A_k \le n_k \le A_k+{\rm M}_\nu,\,\,\, k = 1,2
$$
(this fact follows from (\ref{enny}).
For each ${\bf n}$ the vector
$$
{\bf m} = \mu_* {\bf n} +{\bf m}_{\nu-1} \in {\cal L}_{\nu-1, \mu_*}
$$
together with $ {\bf m}_\nu$ generates a complete lattice
$\langle {\bf m}_\nu, {\bf m}\rangle_\mathbb{Z}$.
Note that $ {\rm M}_\nu \mu_* \asymp H_\nu \cdot \frac{d_{\nu-1}}{M_\nu^\tau M_{\nu-1}} =
o(H_\nu)$ (we use the upper bound from (\ref{volume})). So the set of all vectors ${\bf m}$ constructed is ``dense''  in the range 
$H_\nu \le |\overline{\bf m}|\le 2H_\nu$.
So one can find such a vector with  ${\bf m}\in {\cal D}_\nu$.
That is why we can easy satisfy the condition {\rm (i)} for  the $\nu$-th step.

Similarly,
as
we have many points ${\bf m} \in {\cal D}_\nu$ satisfying {\rm (i)}
we can take ${\bf m}_{\nu+1}$ close enough to $ {\bf w}_\nu$ to satisfy {\rm (iv), (v)}.

As we have certain choice for the vector ${\bf m}_\nu $ at each step of the inductive construction we
can
get $(\alpha_1, \alpha_2) $ 
satisfying linearly independence condition.
So the point $(\alpha_1, \alpha_2)$ constructed satisfies all the conditions of Fundamental Lemma.
The inductive pocedure is described.


\begin{thebibliography}{1}
 \bibitem{SCH}
W.M. Schmidt,\,{\it
Two questions in Diophantine approximations},\,
Monatshefte f\"ur Mathematik\, {\bf 82}, 237 - 245 (1976). 
 \bibitem{SCH1}
  W.M. Schmidt,\, {\it
 Open problems in Diophantine approximations}, in  
 "Approximations Diophantiennes et nombres transcendants`` Luminy, 1982, Progress in
 Mathematics,
 Birkh\"auser, p.271 - 289 (1983).
\bibitem{moshe}
N.G. Moshchevitin,\,
{\it Diophantine approximations with positive integers: a remark to W.M. Schmidt's theorem},
preprint available at arXiv:0904.1906 (2009).


\bibitem{UMN}
N.G. Moshchevitin:\,\,\, Khintchine's singular Diophantine systems and their
applications.// Russian Mathematical Surveys. 65:3 43 - 126 (2010).
\end{thebibliography}
\end{document}